\documentclass[12pt,reqno] {amsart}
\usepackage{graphicx}
\usepackage{tikz}
\usepackage{latexsym, enumitem, fancyhdr, bm}
\usepackage[utf8]{inputenc}
\usepackage{amsfonts}
\usepackage{mathrsfs}
\usepackage{xcolor}
\usepackage{todonotes}
\usepackage{ulem}
\usepackage{hyperref}
\usepackage{mathtools}
\usepackage{amsmath}
\usepackage{amssymb}
\usepackage{amsthm}
\newtheorem{proposition}{Proposition}
\newtheorem{lemma}{Lemma}
\newtheorem{theorem}{Theorem}

\theoremstyle{remark}
\newtheorem{remark}[theorem]{Remark}
\theoremstyle{remark}
\theoremstyle{definition}
\newtheorem{question}{Question}
\newtheorem{definition}[theorem]{Definition}

\title{\textbf{von Neumann inequality for Varopoulos class of operators}}
\author[S. Dey]{Soumyadip Dey}
 \address[S. Dey]{School of Mathematics and Computer Science\\Indian Institute of Technology Goa, India}
 \email{soumyadip22232101@iitgoa.ac.in}

\begin{document}
\begin{abstract}
In this article, we investigate the von Neumann inequality for a class of tuple of commuting bounded linear operators namely Varopoulos operators. In the case of polydisc $\mathbb D^d,$ we show that von Neumann inequality holds for these class of operators with constant $\frac{K^{\mathbb{C}}_G k^{k/2}(k+1)^{(k+1)/2}}{k!\,2^{(k+2)/2}},$ where $K^{\mathbb{C}}_G$ denotes the complex Grothendieck constant and $k$ denotes the degree of the polynomial. For the case of euclidean unit ball $\mathbb B_d,$ we show that von Neumann inequality holds true with constant $1$ for any row contractive tuple of Varopoulos operators. 
\end{abstract}

\maketitle

\section{Introduction}

The classical von Neumann inequality \cite{vN} is one of the fundamental results in operator theory. It asserts that if $T$ is a contraction on a complex Hilbert space and $p\in\mathbb{C}[z]$, then
$\|p(T)\|\leq\|p\|_{\mathbb{D},\infty}$,
where $\|p\|_{\Omega,\infty}$ denotes the supremum of $|p|$ over a domain $\Omega$. Whenever the domain is clear from the context, we simply write $\|p\|_{\infty}$ in place of $\|p\|_{\Omega,\infty}$. A dilation-theoretic proof of this inequality follows from the Sz.-Nagy dilation theorem \cite{Sz-dilation}, which states that every contraction on a Hilbert space admits a unitary dilation. Consequently, $p(T)$ can be realized as a compression of $p(U)$ to the original Hilbert space, yielding the desired norm estimate.

The natural multivariable analogue of the von Neumann inequality has led to a rich and active area of research. In several variables, two natural domains that arise are the unit polydisc $\mathbb{D}^d$ and the Euclidean unit ball $\mathbb{B}_d$.

For the case of polydisc $\mathbb D^d,$ the following is a very natural formulation of von Neumann inequality:
For a commuting $d$-tuple of contractions $\boldsymbol{T}=(T_1,\ldots,T_d)$, whether
$\|p(\boldsymbol T)\|\leq C\|p\|_{\mathbb{D}^d,\infty}$
holds for every polynomial $p\in\mathbb{C}[z_1,\ldots,z_d]$. The case $d=2$ has a particularly satisfactory answer. Indeed, And\^{o}'s dilation theorem \cite{Ando} implies that every pair of commuting contractions $(T_1,T_2)$ admits a commuting unitary dilation and, consequently,
$\|p(T_1,T_2)\|\leq\|p\|_{\mathbb{D}^2,\infty}$.
However, this phenomenon does not persist in higher dimensions.

For the case of euclidean unit ball: A ball analogue of the von Neumann inequality can have various formulations. For instance, one might consider the von Neumann inequality for the class of row contractions or for the class of spherical contractions. However, for both of these classes, the von Neumann inequality does not hold, as evidenced by the works of  Drury (1978)\cite{Drury1978} and Hartz (2022)\cite{Hartz2022}. 
 In \cite{DeyGuptaKumar2026}, authors considered von Neumann inequality for doubly contractive tuples, that is, tuples that are simultaneously row contractions and spherical contractions. They obtained an affirmative result for a subclass of the class of weighted shifts. However, the concerned inequality is not true, in general, even for doubly contractive tuples (see \cite[Theorem 2]{V1}). The main objective of this article is to study von Neumann inequalities, in both polydisc setting and euclidean unit ball setting, for the class of Varopoulos operators (see \cite{Sarkar}).

Let $\mathscr{C}_d$ denote the class of all commuting $d$-tuples of contractions acting on a Hilbert space. Parrott \cite{Parrott} constructed a triple of commuting contractions which does not admit a commuting unitary dilation, although the von Neumann inequality still holds for this particular tuple. Subsequently, Varopoulos \cite{V1}, and independently Crabb and Davie \cite{CD}, showed that the multivariable von Neumann inequality itself fails in general. More precisely, they constructed three commuting contractions $T_1,T_2,T_3$ and a polynomial $p$ such that
$\|p(T_1,T_2,T_3)\|>\|p\|_{\mathbb{D}^3,\infty}$.
Thus, neither the validity of the von Neumann inequality with constant $1$ nor the existence of a commuting unitary dilation can be expected for arbitrary commuting operator tuples in three or more variables.

The failure of the classical von Neumann inequality in several variables naturally leads to the question of whether one can obtain suitable subclasses of commuting contractions for which a von Neumann-type of inequality is true. In this direction, multilinear forms and tensor products provide powerful tools for obtaining  such subclasses for commuting contraction. Of particular relevance are the multidimensional versions of Grothendieck's inequality, in the polydisc settings, and ball analog of Grothendieck inequality developed by Tonge (see \cite[A(r) Inequality]{TG} and \cite[Lemma 4.1]{ManteroTonge1976}). These inequalities provide estimates for multilinear forms which are independent of the number of variables and have proven useful in the study of von Neumann-type inequalities for  classes of Varopoulos tuples (see \cite{Sarkar}).
In \cite[Theorems 1 and 2]{V1}, Varopoulos showed that the von Neumann
inequality cannot, in general, be extended to arbitrary commuting
contractions as well as doubly contractive , when the degree of the homogeneous polynomial is $3$ and
the number of commuting contractions is allowed to increase. 

\begin{question}\label{On-polydisc}
Let $p(z_1,\ldots,z_d)$ be a homogeneous polynomial of degree $k$, and let
$V_1,\ldots,V_d$ be contractive Varopoulos operators of order $k$. Does there exist a
constant $C(k)$, depending only on $k$ and independent of the number of
variables $d$, such that
\[
\left\|p(V_1,\ldots,V_d)\right\|
\leq C(k)\|p\|_{\mathbb{D}^d,\infty}?
\]
\end{question}
\begin{question}\label{Ball-version-of-Varopoulos-Constant-order-k}
Let $p(z_1,\ldots,z_d)$ be a homogeneous polynomial of degree $k$, and let
$V_1,\ldots,V_d$ be row contractive Varopoulos operators of order $k$.
Does it follow that
$\|p(V_1,\ldots,V_d)\|\leq\|p\|_{\mathbb{B}_{d},\infty}$?
\end{question}

The main aim of this article is to investigate the questions stated in Questions \ref{On-polydisc} and \ref{Ball-version-of-Varopoulos-Constant-order-k} (see Theorem \ref{Varopoulos-Constant-order-k}, Theorem \ref{ Ball-version-of-Varopoulos-Constant-order-k}). To this end, we develop several tools that combine techniques from several complex variables, multilinear analysis, and operator theory. In particular, we make use of multidimensional versions of Grothendieck's inequality (see \cite[A(r) Inequality]{TG} and Theorem \ref{Ball-version-of-multidimentional-Grothendieck-inequality}), higher-order derivative estimates for holomorphic functions ( see Theorem \ref{Injective-tensor-norm}), and related estimates for homogeneous polynomials. \begin{definition}[Grothendieck Constant]
For a complex (real) array
$A=(a_{jk})_{d\times d}$, define the norm
\[
\|A\|_{\infty\to 1}
:=
\sup\left\{
|\langle Av,w\rangle|:
\|v\|_{\ell^\infty(d)}\leq 1,\,
\|w\|_{\ell^\infty(d)}\leq 1
\right\},
\]
where $v,w\in\mathbb{C}^d$ (respectively, $\mathbb{R}^d$).
There exists a finite constant $K>0$ such that, for any choice of
unit vectors $(x_j)_{j=1}^{d}$ and $(y_k)_{k=1}^{d}$ in a complex
(respectively, real) Hilbert space $\mathcal{H}$, we have
\[
\left|
\sum_{j,k=1}^{d}a_{jk}\langle x_j,y_k\rangle
\right|
\leq
K\|A\|_{\infty\to 1}
\]
for all $d\in\mathbb{N}$ and $A=(a_{jk})$.
The least such constant is denoted by $K^{\mathbb{C}}_G$ and is known as the complex Grothendieck constant.
\end{definition}

In our context, we use the following  definition to estimate higher-order derivatives and obtain the desired bounds.
\begin{definition}(See \cite{Ryan2002})
Let $X$ and $Y$ be Banach spaces. For an element
$u=\sum_{i=1}^{d}x_i\otimes y_i\in X\otimes Y$, the \textit{injective tensor norm}, denoted by $\|\cdot\|_{\varepsilon}$, is defined by
\[
\|u\|_{\varepsilon}
=\sup_{\substack{x^*\in X^*,\,y^*\in Y^*\\
\|x^*\|\leq 1,\,\|y^*\|\leq 1}}
\Big|\sum_{i=1}^{d}x^*(x_i)y^*(y_i)\Big|.
\]
Equivalently, associate with $u$ the bounded linear operator
$T_u:X^*\to Y$ defined by
\[
T_u(x^*)=\sum_{i=1}^{d}x^*(x_i)y_i,
\qquad x^*\in X^*.
\]
Then
\[
\|u\|_{\varepsilon}=\|T_u\|
=\sup_{\substack{x^*\in X^*\\\|x^*\|\leq 1}}
\Big\|\sum_{i=1}^{d}x^*(x_i)y_i\Big\|_Y.
\]
The completion of $X\otimes Y$ with respect to $\|\cdot\|_{\varepsilon}$ is denoted by $X\widehat{\otimes}_{\varepsilon}Y$.
\end{definition}

The following result, concerning polynomials of degree at most $2$,
provides a motivation for studying von Neumann-type inequalities for
the Varopoulos class of operators and homogeneous polynomials of
arbitrary degree. In particular, it gives a dimension-independent
estimate in terms of the supremum norm on the polydisc.

\begin{theorem}\cite[Theorem 2.1]{RG}
Suppose that $p$ is a polynomial of degree at most $2$ in $n$ variables
and $\boldsymbol{T}\in\mathscr{C}_d$. Then
\[
\|p(\boldsymbol{T})\|
\leq
\frac{3\sqrt{3}}{4}K^{\mathbb{C}}_G
\|p\|_{\mathbb{D}^d,\infty}.
\]
\end{theorem}

\section{Generalized Grothendieck Inequality and Higher-Order Derivative Estimates}
The following result extends the multidimensional Grothendieck inequality of order $3$, we believe the proof presented here has simpler language and  essentially same as that in \cite[A(r) inequality]{TG}. For any Banach space $X,$ let $(X)_1$ denote the open unit ball in that Banach space.
\begin{definition}\label{Symmetric-product}
Let $k$ be positive integers, and let
$x^{(r)}\in(\ell^2)_{1}$, $1\leq r\leq k$.
We define symmetric product of $x^{(1)},x^{(2)},\ldots,x^{(k)}$ as follows
\[
[x^{(1)},\ldots,x^{(k)}]
:=
\sum_{l=1}^{\infty}
x_{l}^{(1)}
\cdots x_{l}^{(k)}.
\]
\end{definition}

\begin{theorem}
Let $d$ be a positive integer, and let $x_{1},\ldots,x_{d},$ $y_{1},\ldots,y_{d},$ $z_{1},\ldots,z_{d} \in (\ell^{2})_1$.
Then for every tensor $b=(b_{ijk})_{i,j,k=1}^{d}$, the following estimate holds:
\[
\left| 
   \sum_{i,j,k=1}^{d} b_{ijk}\,[x_i,y_j,z_k]
\right|
   \leq K^{\mathbb{C}}_G\,\sqrt{2}\;
      \| b \|_{\,\ell^{1}\otimes_{\varepsilon} \ell^{1}\otimes_{\varepsilon} \ell^{1}} .
\]
\end{theorem}
\begin{proof}
Define the operator $U : \ell^\infty(d) \longrightarrow \ell^1(d)$
by prescribing its action on basis vectors:
$$e_m \longmapsto u_m:=\sum_{i,j=1}^d b_{ijm}\,(y_i \circ z_j),$$
where $(y_i\circ z_j)^{(n)} = y_i^{(n)} z_j^{(n)}$. 
In these notations, note that
\begin{eqnarray}\label{Littlewood's-Application}
   \big|\sum_{i,j,m=1}^d b_{ijm}[x_m,y_i,z_j]\big| = \big|\sum_{m=1}^d \langle x_m , u_m\rangle\big|\leq \sum_{m=1}^{d}\big(\sum_{n=1}^{\infty}|u_m^{(n)}|^2\big)^{1/2}\le \sqrt{2}\|U\|_{\ell^\infty(d)\to \ell^1(d)},
\end{eqnarray}
where the last inequality follows from Littlewood's inequality (see \cite[page-281]{Li} and \cite{Sz}).
Further, we note that 
\begin{eqnarray*}
    \|U\|
   &=& \sup_{\|\alpha\|_\infty\le 1, \|\beta\|_\infty\le 1}
      \Big|
        \sum_{n,m} u_m^{(n) }\alpha_n \beta_m
      \Big|\\ &=& \sup_{|\alpha_n|=1,\;|\beta_m|=1}
   \left|
     \sum_{i,j,m}
        b_{ijm}\,\beta_m
     \sum_n y_i^{(n)} z_j^{(n)} \alpha_n
   \right|\\
 &=& \sup_{|\alpha_n|=1,\;|\beta_m|=1}
   \left|
     \sum_{i,j} a_{ij}\,
        \langle \tilde y_i , z_j\rangle
   \right|,
\end{eqnarray*}
  where for each $i,j=1,\ldots, d,$ $a_{ij}=\sum_{m=1}^d \beta_m\,b_{ijm}$ and $\tilde y_i = \alpha \circ y_i$. 
  From Grothendieck's inequality, we get that 
  $$\|U\|\leq K_G^{\mathbb C} \|(\!(a_{ij})\!)\|_{\ell^1(d)\otimes_\varepsilon\ell^1(d)}=K_G^{\mathbb C}\sup_{|\chi_i|=|\eta_j|=1}\big|\sum_{i,j} a_{ij}\,\chi_i \eta_j\big| \leq K_G^{\mathbb{C}}\sup_{|\chi_i|=|\eta_j|=|\beta_m|=1}\Big|\sum_{i,j,m=1}^{d}b_{ijm}\chi_i \eta_j \beta_m\Big|.$$ 
  Using this in \eqref{Littlewood's-Application}, we get,
  \begin{eqnarray*}
      \Big| \sum_{i,j,m} b_{ijm}[x_m,y_i,z_j] \Big|
 \leq
K^{\mathbb{C}}_G\sqrt{2} \|b\|_{\,\ell^1\otimes_\varepsilon\ell^1\otimes_\varepsilon\ell^1}.
  \end{eqnarray*}
This completes the proof.
\end{proof}
The following estimate for higher-order derivatives of bounded holomorphic functions will be useful in our analysis. Although this estimate is known from the work of Harris (see \cite{Harris1975}), we provide an alternative proof based on the Schwarz lemma. The proof presented here is along the idea presented in \cite{Gupta2017}. 
\begin{theorem}\label{Injective-tensor-norm} 
Let  $f: \mathbb{D}^{d} \to \mathbb{D}$ be a holomorphic function. Then, for the $k$th derivative of $f$ at $0$, we have the following bound:
\[
 \|D^{k} f(0)\|_{\ell^1(d){\otimes}_\varepsilon \ell^1(d){\otimes}_\varepsilon \cdots {\otimes}_\varepsilon \ell^1(d)} \leq \frac{k^{k/2}(k+1)^{(k+1)/2}}{2^{k}},
\]
where $ k \geq 1 $.
\end{theorem}
\begin{proof}
The case $k=1$ follows directly from the Schwarz's lemma (see \cite[Theorem 8.1.2]{Rudin}). The case $k=2$ is done in \cite{Gupta2017}.
 Let $k=3$.
Let $f$ be a complex-valued analytic function on $ \mathbb{D}^d $ with $ \lVert f \rVert_{\mathbb{D}^n, \infty} \leq 1.$  
Consider an arbitrary point $ a = (a_1, \dots, a_d) \in r_1 \mathbb{D}^d $, where $ 0 < r_1 < 1 .$ Let $ \Phi_j $ be an automorphism of $ r_1 \mathbb{D} $, defined by $\Phi_j(z) = r_1^2 \frac{a_j - z}{r_1^2 - \overline{a_j}z}$ for  $j = 1, \dots, n.$
Define $ \Phi(z_1, \dots, z_d) := (\Phi_1(z_1), \dots, \Phi_d(z_d)) $, an automorphism on $ r_1 \mathbb{D}^d $, which maps $0$ to $a$. Since $ f $ is holomorphic, $ Df $ is also holomorphic, mapping $ r_1 \mathbb{D}^d $ to $ \frac{1}{1 - r_1^2} (\mathbb{D}^d)^* $, where $(\mathbb D^d)^*$ denotes the open unit ball of $\ell^1(d)$ (dual of $\ell^\infty(d)$). 
Define the holomorphic function $$g := Df \circ \Phi : r_1 \mathbb{D}^d \to \frac{1}{1 - r_1^2} (\mathbb{D}^d)^*.$$ Then, $ Dg(0) = D^2f(a) D\Phi(0) $, which in turn gives $D^2f(a) = Dg(0) \left( D\Phi(0) \right)^{-1}.$ 
Thus, we have the bound: $\lVert D^2f(a) \rVert_{\infty \to 1} \leq \lVert Dg(0) \rVert_{\infty \to 1} \lVert (D\Phi(0))^{-1} \rVert_{\infty \to \infty}.$
By Schwarz's lemma, $ Dg(0) $ is a linear transformation mapping $ r_1 \mathbb{D}^d $ to $ \frac{1}{1 - r_1^2} (\mathbb{D}^d)^* $, therefore
\[
\lVert Dg(0) \rVert_{\infty \to 1} \leq \frac{1}{r_1(1 - r_1^2)}.
\]
As \( D\Phi(0) = \operatorname{diag} \left( \frac{\lvert a_1 \rvert^2 - r_1^2}{r_1^2}, \dots, \frac{\lvert a_d \rvert^2 - r_1^2}{r_1^2} \right) \), we get
\[
\lVert (D\Phi(0))^{-1} \rVert_{\infty \to \infty} \leq \max_j \frac{r_1^2}{r_1^2 - \lvert a_j \rvert^2}.
\]
Let $ 0 < r_2 < 1 $ such that $ r_1r_2 = \max_j \lvert a_j \rvert $. Then, $ \lvert a_j \rvert \leq r_1r_2 $, and hence
\[
\lVert (D\Phi(0))^{-1} \rVert_{\infty \to \infty} \leq \frac{1}{1 - r_2^2}.
\]
Therefore,
\[
\lVert D^2f(a) \rVert_{\infty \to 1} \leq \frac{1}{r_1(1 - r_1^2)(1 - r_2^2)}.
\]
Thus, $ D^2f : r_1r_2 \mathbb{D}^d \to \frac{1}{r_1(1 - r_1^2)(1 - r_2^2)} (\ell^1(d) {\otimes}_\varepsilon \ell^1(d))_1$ is holomorphic. By Schwarz's Lemma,
$D^3f(0)$ is a linear map that maps  $r_1r_2 \mathbb{D}^d$ into $\frac{1}{r_1(1 - r_1^2)(1 - r_2^2)} (\ell^1(d) {\otimes}_\varepsilon \ell^1(d))_1.$
Hence,
\[
\lVert D^3f(0) \rVert_{\ell^\infty(d) \to (\ell^1(d) {\otimes}_\varepsilon \ell^1(d))} \leq \frac{1}{r_1^2(1 - r_1^2)r_2(1 - r_2^2)}.
\]
Taking infimum over $r_1$ and $r_2$ on the right hand side, we get
$\lVert D^3f(0) \rVert_{\ell^\infty(d) \to (\ell^1(d) {\otimes}_\varepsilon \ell^1(d))} \leq 6 \sqrt{3}.$ This completes the proof of theorem for $k=3.$

  For general $k\in\mathbb N,$ we claim the following first:

  \textbf{Claim:} Suppose $ f : \mathbb{D}^d \to \mathbb{D} $ is a holomorphic map. Then, for any fixed $ k \in \mathbb{N} $ and $ r_1, \dots, r_k \in (0,1) $,
\[
D^kf : r_1r_2 \cdots r_k \mathbb{D}^d \to \alpha_k (\ell^1(d) {\otimes}_\varepsilon \cdots {\otimes}_\varepsilon \ell^1(d))_1 \quad (k\text{-times}),
\]
where
\[
\alpha_k = \frac{1}{r_1^{k-1}(1 - r_1^2) \cdots r_{k-1}(1 - r_{k-1}^2)(1 - r_k^2)}.\]
\textbf{Proof of the claim:} 
The proof proceeds by induction on $ k $. For $ k = 1 $, the result holds trivially since $ Df $ is holomorphic on $ r_1 \mathbb{D}^d $. Now, assume the is true for some $ k \geq 1 $. We show that it holds for $ k+1 $. 
Consider the composition:
\[
r_1r_2 \cdots r_k \mathbb{D}^d \stackrel{\Phi}{\longrightarrow} r_1r_2 \cdots r_k \mathbb{D}^d \stackrel{D^kf}{\longrightarrow} \alpha_k (\ell^1(d) {\otimes}_\varepsilon \cdots {\otimes}_\varepsilon \ell^1(d))_1.
\]
where $ \Phi(z_1, \dots, z_d) = (\Phi_1(z_1), \dots, \Phi_d(z_d)) $, an automorphism on $ r_1r_2\cdots r_k \mathbb{D}^d $, of the form \[
\Phi_j(z_j)=R^2\frac{a_j-z_j}{R^2-\overline{a_j}z_j},
\qquad R=r_1r_2\cdots r_k,
\] which maps $0$ to $a$. 
Let $ g = D^kf \circ \Phi $. Then, $ g $ is holomorphic, and $ Dg(0) = D^{k+1}f(a) D\Phi(0) $. Therefore, $ D^{k+1}f(a) = Dg(0) (D\Phi(0))^{-1} $.  By Schwarz's lemma, we have that the differential $ Dg(0) $ maps
 $r_1 r_2 \cdots r_k \mathbb{D}^d$ into $\alpha_k \left( (\ell^1(d) {\otimes}_\varepsilon \ell^1(d) {\otimes}_\varepsilon \cdots {\otimes}_\varepsilon \ell^1(d))_1 \right).$
 It follows that
\[ \| Dg(0) \|_{\ell^{\infty}(d) \to (\ell^1(d) {\otimes}_\varepsilon \ell^1(d) {\otimes}_\varepsilon \cdots {\otimes}_\varepsilon \ell^1(d))} \leq \frac{1}{r_1^k (1 - r_1^2) \cdots r_k (1 - r_k^2)}. \]
As $ \Phi $ is holomorphic and $ (D \Phi(0))^{-1} $ is invertible, then
\[ \| (D \Phi(0))^{-1} \|_{\infty \to \infty} \leq \frac{1}{1 - r_{k+1}^2}, \]
where \( 0 < r_{k+1} < 1 \) is such that \( r_1 \cdots r_{k+1} = \max_j |a_j| \).
Thus,
\[ D^{k+1} f : r_1 \cdots r_{k+1} \mathbb{D}^d \rightarrow \frac{1}{r_1^k (1 - r_1^2) \cdots r_k (1 - r_k^2) (1 - r_{k+1}^2)} (\ell^1(d) {\otimes}_\varepsilon \ell^1(d) {\otimes}_\varepsilon \cdots {\otimes}_\varepsilon \ell^1(d))_{1}, \]
and the claim stands verified.

Since $ D^k f $ is holomorphic, $ D^{k+1} f $ is also holomorphic. Applying Schwarz's lemma, we get
\[ \| D^{k+1} f(0) \|_{\ell^{\infty}(d) \to (\ell^1(d) {\otimes}_\varepsilon \ell^1(d) {\otimes}_\varepsilon \cdots {\otimes}_\varepsilon \ell^1(d))} \leq \frac{1}{r_1^k (1 - r_1^2) \cdots r_k (1 - r_k^2)}. \]

To minimize the right-hand side, we maximize the denominator of the fraction. First, we calculate
\[ \max r^k (1 - r^2), \]
which, by elementary calculus, is maximized at $ r = \sqrt{\frac{k}{k+2}} $. Thus,
\[ \| D^{k+1} f(0) \|_{\ell^{\infty}(d) \to (\ell^1(d) {\otimes}_\varepsilon \ell^1(d) {\otimes}_\varepsilon\cdots {\otimes}_\varepsilon \ell^1(d))} \leq \frac{(k+1)^{(k+1)/2} (k+2)^{(k+2)/2}}{2^{k+1}}. \]
This completes the proof of the theorem.
\end{proof}
\begin{definition}\cite{P.G Dixon}
$BM(k)$-norm of a symmetric $k$-tensor $a=(a_{i_1,\ldots,i_k})$ is given by
\[
\|a\|_{BM(k)}:=\sup\Bigg\{
\Big|
\sum_{i_1,\ldots,i_k=1}^d a_{i_1,\dots,i_k}\, f_1(i_1)\cdots f_k(i_k)
\Big|
:\ 
f_i\in C(I_d),\ \|f_i\|\leq 1,\ 1\leq i\leq k 
\Bigg\}
\]
where $I_d=\{1,\ldots , d\}.$
\end{definition}
Corresponding to any symmetric $k$-tensor $a=(a_{i_1,\ldots,i_k}),$ one can associate a homogeneous polynomial $p$ of degree $k$ in $d$-variables defined by
\begin{eqnarray}\label{polynomial-corresponding-to-tensor}
    p(z_1,\dots,z_d)
   = \sum_{i_1,\ldots,i_k=1}^d a_{i_1,\ldots,i_k}\, z_{i_1}\cdots z_{i_k}.
\end{eqnarray}
Clearly $\|p\|_{\mathbb D^d,\infty} \le \|a\|_{BM(k)}$.
The following lemma gives the reverse bound and can be found in \cite{DV}.
\begin{lemma}\cite[Lemma 2.1]{DV}
If $a$ is a symmetric tensor of order $k$ determining the polynomial $p$, then
\begin{eqnarray}\label{BM-1}
\|p\|_\infty
   \;\leq\;
\|a\|_{BM(k)}
   \;\leq\;\frac{(2^k-1)\,k^k}{k!}\,\|p\|_\infty
\;\leq\;(2e)^k\,\|p\|_\infty.
\end{eqnarray}
\end{lemma}
Theorem \ref{Injective-tensor-norm} improves the right hand side bound in \eqref{BM-1} considerably which is noted down in the following lemma.
\begin{lemma}
If $a$ is a symmetric tensor of order $k$ determining the polynomial $p$ (as given in \eqref{polynomial-corresponding-to-tensor}), then \begin{eqnarray}\label{BM-2}
\|p\|_\infty
   \;\leq\;
\|a\|_{BM(k)}
   \;\leq\;\frac{k^{k/2} (k+1)^{(k+1)/2}}{2^kk!}\|p\|_\infty.
\end{eqnarray}
\end{lemma}
The following remark compares the constants appearing on the right-hand sides of \eqref{BM-1} and \eqref{BM-2} and shows that \eqref{BM-2} gives a much sharper bound.
\begin{remark}
Let $k\in\mathbb N$ and
\[
A_k=\frac{(2^k-1)k^k}{k!},
\qquad
B_k=\frac{k^{k/2}(k+1)^{(k+1)/2}}{2^k k!}.
\]
By Stirling's formula, one gets 
\[
A_k\sim \frac{(2e)^k}{\sqrt{2\pi k}}.
\]
Moreover, using $(1+1/k)^{k/2}\to e^{1/2}$, we obtain
\[
\begin{aligned}
B_k
&=\frac{k^{k+1/2}}{2^k k!}
  \left(1+\frac1k\right)^{(k+1)/2}\sim \frac{e^k}{2^k\sqrt{2\pi}}.
\end{aligned}
\]
It follows that
\[
\frac{A_k}{B_k}\sim \frac{4^k}{\sqrt{k}},
\]
which tends to infinity as $k\to\infty$. Hence the bound given by $B_k$ is much better than that by $A_k$ for all sufficiently large $k$.
\end{remark}

The following statement is motivated by the claim proved in Theorem \ref{Injective-tensor-norm}. It shall be used in Theorem \ref{Gneneralize-injective-norm}.
\begin{proposition}\label{Bounds-on-ball}
Suppose $ f : \big(\ell^{\infty}(d)\otimes_{\varepsilon}\ell^2(d)\big)_1 \to \mathbb{D} $ is a holomorphic map. Then, for any fixed $ k \in \mathbb{N} $ and $ r_1, \dots, r_k \in (0,1) $,
\[
D^kf : r_1r_2 \cdots r_k (\ell^{\infty}(d)\otimes_{\varepsilon}\ell^2(d))_1 \to \alpha_k (\mathcal{L}_{k})_1,
\]
where
\[
\alpha_k = \frac{1}{r_1^{k-1}(1 - r_1^2) \cdots r_{k-1}(1 - r_{k-1}^2)(1 - r_k^2)},
\] and
\[ \mathcal{L}_{k}=\mathcal{L}(\ell^{\infty}(d)\otimes_{\varepsilon}\ell^2(d),\mathcal{L}_{k-1})\] 
with $\mathcal{L}_{0}=\mathbb{C}$.
\end{proposition}
In order to prove Proposition \ref{Bounds-on-ball}, we need some technique.
We work in the Banach space
$\ell^{\infty}(d)\otimes_{\varepsilon}\ell^{2}(d).$
Let $x_1,\ldots,x_d$ be vectors in $(\ell^{2}(d))_1$, and consider the element $x=(x_1,\ldots,x_d)\in
(\ell^{\infty}(d)\otimes_{\varepsilon}\ell^{2}(d))_1.$ Note that the injective norm of $x$ in this setting is given by 
\[
\|x\|_{\ell^{\infty}(d)\otimes_{\varepsilon}\ell^{2}(d)}
   = \max_{1\leq i\leq d}\|x_i\|_{\ell^{2}(d)}.
\]
Let 
$f : \bigl(\ell^{\infty}(d)\otimes_{\varepsilon}\ell^{2}(d)\bigr)_1\to\mathbb{C}$
be analytic with $\|f\|_{\infty}\le 1.$
Fix $a=(a_1,\dots,a_d)$ with each $a_i\in (\ell^{2}(d))_1$.
Define the automorphism $\Phi$ on $\bigl(\ell^{\infty}(d)\otimes_{\varepsilon}\ell^{2}(d)\bigr)_1$ given by
\begin{eqnarray}\label{Automorphism-form-on-ball}
\Phi(x_1,\dots,x_d)
  = \bigl(\phi_{a_1}(x_1),\dots,\phi_{a_d}(x_d)\bigr),
\end{eqnarray}
where each $\phi_{a_i}$ is the usual euclidean unit ball automorphism defined by, for any $y\in(\ell^2)_1,$
\[
\phi_{a_i}(y)
   = \frac{a_i - P_{a_i}(y) - s_{a_i} Q_{a_i}(y)}
          {1 - \langle y, a_i\rangle},
\]
with
\[
P_{a_i}(y)
   := \frac{\langle y,a_i\rangle}{\langle a_i,a_i\rangle}a_i,
\qquad
Q_{a_i}(y)
   := y - \frac{\langle y,a_i\rangle}{\langle a_i,a_i\rangle}a_i,
\qquad
s_{a_i} := \sqrt{1 - \|a_i\|^2}.
\]
Set $g=f\circ \Phi$, which is analytic on the unit ball. Note that 
$Dg(0)=Df(\Phi(0))\circ D\Phi(0),$ by the chain rule. 
The derivative of $\Phi$ at $0,$ $D\Phi(0),$ is block diagonal $D\Phi(0)
  = \operatorname{diag}(D\phi_{a_1}(0),\dots,D\phi_{a_d}(0)).$
For each $i\in\{1,\ldots, d\},$
$D\phi_{a_i}(0)
   = -s_{a_i}^{2} P_{a_i} - s_{a_i} Q_{a_i}$ (see \cite[Page 25]{Rudin}).
Note that $M_j := D\phi_{a_j}(0)$ is invertible.  Let $U_j$ denote the unitary matrix which diagonalizes $P_{a_j}$ and $Q_{a_j}$ simultaneously. Then,
\[
U_j M_j U_j^{*}
   = \begin{pmatrix}
       -s_{a_j}^{2} & 0 & \cdots & 0 \\
       0 & -s_{a_j} & \cdots & 0 \\
       \vdots & \vdots & \ddots & \vdots \\
       0 & 0 & \cdots & -s_{a_j}
     \end{pmatrix}.
\]
Thus,
\[
M_j^{-1}
=
U_j^{*}
\begin{pmatrix}
-\dfrac{1}{s_{a_j}^{2}} & 0 & \cdots & 0 \\
0 & -\dfrac{1}{s_{a_j}} & \cdots & 0 \\
\vdots & \vdots & \ddots & \vdots \\
0 & 0 & \cdots & -\dfrac{1}{s_{a_j}}
\end{pmatrix}
U_j,
\]
and
\[
D\Phi(0)^{-1}
   = \operatorname{Diag}(M_1^{-1},\dots,M_d^{-1}).
\]
By Schwarz's lemma, note that  $Dg(0):(\ell^{\infty}(d)\otimes_{\varepsilon}\ell^2(d))_1 \to \mathbb{D}$  is a holomorphic map. 
Writing this blockwise yields
\[Df(a)=(Dg(0)M_1^{-1},\dots,Dg(0)M_d^{-1}).\]
Therefore the operator norm is bounded by
\[\|Df(a)\|_{(\ell^{\infty}(d)\otimes_{\varepsilon}\ell^{2}(d))_1^{*}}\leq \|Dg(0)\|_{(\ell^{\infty}(d)\otimes_{\varepsilon}\ell^{2}(d))_1^{*}}\cdot \max_{1\leq j\leq d}\Big\{\frac{1}{s_{a_j}^{2}},\frac{1}{s_{a_j}}\Big\}.
\]
If $r_1\in(0,1)$ satisfies $\|a_j\|<r_1$ for all $j$, then $s_{a_j}>\sqrt{1-r_1^{2}}$ and one obtains the estimate
\[\|Df(a)\|_{(\ell^{\infty}(d)\otimes_{\varepsilon}\ell^{2}(d))_1^{*}}\leq\frac{1}{1-r_1^{2}}.\]
 Thus the function $Df$ is also holomorphic that maps $r_1(\ell^{\infty}(d)\otimes_{\varepsilon}\ell^2(d))_1 $ to $\frac{1}{1-r_1^2}(\ell^{\infty}(d)\otimes_{\varepsilon}\ell^2(d))^{*}_1.$
Therefore once again by Schwarz's lemma $D^2f(0)$ is a linear transformation that maps
$r_1(\ell^{\infty}(d)\otimes_{\varepsilon}\ell^2(d))_1 $ into $\frac{1}{1-r_1^2}(\ell^{\infty}(d)\otimes_{\varepsilon}\ell^2(d))^{*}_1$ 
and hence, 

\[\|D^2f(0)\|_{(\ell^{\infty}(d)\otimes_{\varepsilon}\ell^2(d))_1 \to (\ell^{\infty}(d)\otimes_{\varepsilon}\ell^2(d))^{*}_1} \leq  \frac{1}{r_1(1-r_1^2)}.\]
Let $g = Df\circ\Phi : r_1(\ell^{\infty}(d)\otimes_{\varepsilon} \ell^2(d))_1 \to
\frac{1}{1-r_1^2}\bigl(\ell^{\infty}(d)\otimes_{\varepsilon} \ell^2(d)\bigr)_1^*$ be holomorphic, where $\Phi$ is an automorphism on $r_1(\ell^{\infty}(d)\otimes_{\varepsilon} \ell^2(d))_1$ which is of the form \eqref{Automorphism-form-on-ball}.
By Schwarz's lemma, $Dg(0)$ is a linear map that maps $r_1(\ell^{\infty}(d)\otimes_{\varepsilon} \ell^2(d))_1$ into
$\frac{1}{1-r_1^2}\bigl(\ell^{\infty}(d)\otimes_{\varepsilon} \ell^2(d)\bigr)_1^*.$ Then
\[
\begin{aligned}
\|D^2 f(a)\|_{
(\ell^\infty(d)\otimes_\varepsilon \ell^2(d))_1
\to
(\ell^\infty(d)\otimes_\varepsilon \ell^2(d))_1^*
}
&\leq
\|Dg(0)\|_{
(\ell^\infty(d)\otimes_\varepsilon \ell^2(d))_1
\to
(\ell^\infty(d)\otimes_\varepsilon \ell^2(d))_1^*
}
\\
&\quad\times
\|(D\Phi(0))^{-1}\|_{
(\ell^\infty(d)\otimes_\varepsilon \ell^2(d))_1
\to
(\ell^\infty(d)\otimes_\varepsilon \ell^2(d))_1
}.
\end{aligned}
\] 
We also have $\|(D\Phi(0))^{-1}\|_{(\ell^{\infty}(d)\otimes_{\varepsilon} \ell^2(d))_1 \to
(\ell^{\infty}(d)\otimes_{\varepsilon} \ell^2(d))_1}\leq\max_j \frac{r_1^2}{r_1^2-\|a_j\|^2}.$
Let $0<r_2<1$ be such that $r_1 r_2 = \max_{1\le j\le n}\|a_j\|$ (such that  $\|a_j\|\leq r_1 r_2$ for all $j$)  and hence,
\[\|(D\Phi(0))^{-1}\|_{(\ell^{\infty}(d)\otimes_{\varepsilon} \ell^2(d))_1 \to
(\ell^{\infty}(d)\otimes_{\varepsilon} \ell^2(d))_1} \leq  \frac{1}{1-r_2^2}.\] 
Therefore, \[\|D^2 f(a)\|_{\,(\ell^{\infty}(d)\otimes_{\varepsilon} \ell^2(d))_1 \to
(\ell^{\infty}(d)\otimes_{\varepsilon} \ell^2(d))_1^*}
\le \frac{1}{r_1(1-r_1^2)(1-r_2^2)}.\]
Thus, $$D^2f:r_1r_2(\ell^{\infty}(d)\otimes_{\varepsilon} \ell^2(d))_1 \to \frac{1}{r_1(1-r_1^2)(1-r_2^2)} \mathcal{L}((\ell^{\infty}(d)\otimes_{\varepsilon} \ell^2(d))_1 ,\bigl(\ell^{\infty}(d)\otimes_{\varepsilon} \ell^2(d)\bigr)_1^*)$$ 
is holomorphic. 
Once again using Schwarz's lemma, we get  
\[\|D^3f(0)\|_{(\ell^{\infty}(d)\otimes_{\varepsilon} \ell^2(d))_1 \to  \mathcal{L}((\ell^{\infty}(d)\otimes_{\varepsilon} \ell^2(d))_1 ,(\ell^{\infty}(d)\otimes_{\varepsilon} \ell^2(d))_1^*)} \leq \frac{1}{r_1^2(1-r_1^2)r_2(1-r_2^2)}.\]
Minimum of right hand side  occurs at \( r_1 = \frac{1}{\sqrt{2}} \) and \( r_2 = \frac{1}{\sqrt{3}} \), therefore we get 
\[\|D^3f(0)\|_{(\ell^{\infty}(d)\otimes_{\varepsilon} \ell^2(d))_1 \to  \mathcal{L}((\ell^{\infty}(d)\otimes_{\varepsilon} \ell^2(d))_1,(\ell^{\infty}(d)\otimes_{\varepsilon} \ell^2(d))_1^*)} \leq  6\sqrt{3}.\]

\begin{proof}[{Proof of Proposition \ref{Bounds-on-ball}:}]
The proof proceeds by induction on $ k $. For $ k = 1 $, the result holds trivially since $ Df $ is holomorphic on $ r_1 (\ell^{\infty}(d)\otimes_{\varepsilon}\ell^2(d))_1 $. Now, assume the result holds for some $ k \geq 1 $. We show that it holds for $ k+1 $. 
Let $g$ denote the following composition:
$$r_1r_2 \cdots r_k (\ell^{\infty}(d)\otimes_{\varepsilon}\ell^2(d))_1 \stackrel{\Phi}{\longrightarrow} r_1r_2 \cdots r_k (\ell^{\infty}(d)\otimes_{\varepsilon}\ell^2(d))_1 \stackrel{D^kf}{\longrightarrow} \alpha_k (\mathcal{L}_{k})_1,$$
 where $\Phi$ is an automorphism on $r_1r_2 \cdots r_k (\ell^{\infty}(d)\otimes_{\varepsilon}\ell^2(d))_1.$ 
Then, $ g $ is holomorphic, and $ Dg(0) = D^{k+1}f(a) D\Phi(0) $. Therefore, $ D^{k+1}f(a) = Dg(0) (D\Phi(0))^{-1} $.  By Schwarz's lemma, we have that the differential $ Dg(0) $ maps
 $r_1 r_2 \cdots r_k (\ell^{\infty}(d)\otimes_{\varepsilon}\ell^2(d))_1$  into $\alpha_k \left( (\mathcal{L}_{k})_1\right)$ and
 
 \[\| Dg(0) \|_{(\ell^{\infty}(d)\otimes_{\varepsilon}\ell^2(d))_1  \rightarrow  \left( (\mathcal{L}_{k})_1\right)} \leq \frac{1}{r_1^k (1 - r_1^2) \cdots r_k (1 - r_k^2)}.\] 
Additionally,  $ \Phi $ is holomorphic and $ (D \Phi(0))^{-1} $ is invertible, then

 \[\| (D \Phi(0))^{-1} \|_{(\ell^{\infty}(d)\otimes_{\varepsilon} \ell^2(d))_1 \to (\ell^{\infty}(d)\otimes_{\varepsilon} \ell^2(d))_1 } \leq \frac{1}{1 - r_{k+1}^2},\] 
where $ 0 < r_{k+1} < 1 $ is such that $ r_1 \cdots r_{k+1} = \max_j \|a_j\| $.
Thus,
\[D^{k+1} f : r_1 \cdots r_{k+1} (\ell^{\infty}(d)\otimes_{\varepsilon} \ell^2(d))_1  \rightarrow \frac{1}{r_1^k (1 - r_1^2) \cdots r_k (1 - r_k^2) (1 - r_{k+1}^2)} (\mathcal{L}_{k+1})_1.\]
This completes the proof.
\end{proof}
The following result provides a generalization of the preceding estimate
for higher-order derivatives (see Theorem \ref{Injective-tensor-norm}).
\begin{theorem}\label{Gneneralize-injective-norm}
  Suppose $ f : \big(\ell^{\infty}(d)\otimes_{\varepsilon}\ell^2(d)\big)_1 \to \mathbb{D} $ is a holomorphic map. Then, for any fixed $ k \in \mathbb{N} $ and $ r_1, \dots, r_k \in (0,1) $, then \[ \| D^k f(0) \|_{(\ell^{\infty}(d)\otimes_{\varepsilon} \ell^2(d))_1 \to (\mathcal{L}_{k-1})_1} \leq \frac{k^{k/2} (k+1)^{(k+1)/2}}{2^k}. \]  
\end{theorem}
\begin{proof}
Using Proposition \ref{Bounds-on-ball}, we have
\[ r_1 r_2 \cdots r_k (\ell^{\infty}(d)\otimes_{\varepsilon} \ell^2(d))_1 \stackrel{D^k f}{\longrightarrow} \frac{1}{r_1^{k-1} (1 - r_1^2) \cdots r_{k-1} (1 - r_{k-1}^2) (1 - r_k^2)} (\mathcal{L}_{k})_1. \]
Since $ D^k f $ is holomorphic, $D^{k+1} f$ is also holomorphic. Applying Schwarz's lemma, we get
\[ \| D^{k+1} f(0) \|_{(\ell^{\infty}(d)\otimes_{\varepsilon} \ell^2(d))_1 \to (\mathcal{L}_{k})_1} \leq \frac{1}{r_1^k (1 - r_1^2) \cdots r_k (1 - r_k^2)}. \]
To minimize the right-hand side, we maximize the denominator of the fraction. First, we calculate $\max r^k (1 - r^2),$
which by elementary calculus, is maximized at \( r = \sqrt{\frac{k}{k+2}} \). Thus,
\[ \| D^{k+1} f(0) \|_{(\ell^{\infty}(d)\otimes_{\varepsilon} \ell^2(d))_1 \to (\mathcal{L}_{k})_1} \leq \frac{(k+1)^{(k+1)/2} (k+2)^{(k+2)/2}}{2^{k+1}}. \]
This completes the proof of the theorem.
\end{proof}
\subsection{Asymptotic Behavior of the Ratio for Commuting Contractions and Varopoulos Operators of order $2$}

Let $f : (\ell^{\infty}(d)\otimes_{\varepsilon}\ell^2(d))_1 \to \mathbb{D}$ be the holomorphic function defined by 
$$f(x_1,\dots,x_d)=\sum_{i,j=1}^{d} a_{ij} [x_i,x_j].$$ 
Then $D^{2}f(0)=2\,(A\otimes I),$ where $I$ is an identity operator on $\ell^2(d)$ and $A=(a_{ij})$.
Hence
\begin{eqnarray*}  \|D^{2}f(0)\|_{(\ell^{\infty}(d)\otimes_{\varepsilon}\ell^{2}(d))_{1}\to (\mathcal{L}_{1})_1}=2\sup_{\substack{\|x\|_{(\ell^{\infty}(d)\otimes_{\varepsilon}\ell^{2}(d))}\leq 1\\\|y\|_{(\ell^{\infty}(d)\otimes_{\varepsilon}\ell^{2}(d))}\leq 1}}\Big|\sum_{i,j=1}^{d} a_{ij} [x_i,y_j]\Big|.
\end{eqnarray*}
Combining with the Theorem \ref{Gneneralize-injective-norm}, we have for all $x=(x_1,\ldots ,x_d)\in (\ell^{\infty}(d)\otimes_{\varepsilon}\ell^{2}(d))_1,$ $y=(y_1,\ldots ,y_d)\in (\ell^{\infty}(d)\otimes_{\varepsilon}\ell^{2}(d))_1$ 
\begin{eqnarray}\label{ratio-compare}
\Big|\sum_{i,j=1}^{d} a_{ij} [x_i,y_j]\Big|\le \frac{3\sqrt{3}}{4}\,\|f\|_{\infty}.
\end{eqnarray}
Let $\mathbf{T}=(T_1,\dots,T_d)$ be a tuple of commuting contractions on a Hilbert space, and let $p$ be a non-zero homogeneous polynomial of degree $2$ in $d$ variables. Then, with some efforts, \eqref{ratio-compare} can be written as
\[\frac{\sup_{\boldsymbol{T}\in \mathscr{C}_d}\|p(T_1,\dots,T_d)\|}{\sup_{\boldsymbol V\in\mathcal{V}_{d,2}}\|p(V_1,\dots,V_d)\|}\leq \frac{3\sqrt{3}}{4},\]
where $$V_i=\begin{pmatrix}
0 & x_i & 0 &  \\[6pt]
0 & 0 & x_i^t  \\[6pt]
0 & 0 & 0 &  \\[6pt]
\end{pmatrix}.$$  Note that $(V_1,\dots,V_d)$ is a commuting tuple of bounded operators known as Varopoulos operators of order $2$ (see \cite{Gupta}) and the set of all commuting contractive $d$ tuple of Varopoulos operators of order $2$ is denoted by $\mathcal{V}_{d,2}.$ 

The following definition of  Varopoulos operators of order $k$ has been taken from \cite[page 189]{Sarkar} 

\begin{definition}[ Varopoulos Operators of order $k$]
Fix $d,k\in\mathbb{N}$ with $k\geq 2$. For any row vectors
$x_1,\ldots,x_d\in\mathbb{C}^d$, define the operators,
$V_j=
\begin{pmatrix}
0 & x_j & 0 & 0\\
0 & 0 & D_{x_j}\otimes I_{k-2} & 0\\
0 & 0 & 0 & x_j^{\mathrm{t}}\\
0 & 0 & 0 & 0
\end{pmatrix},$ for $j=1,\ldots,d,$
where $D_{x_j}$ denotes the diagonal matrix with diagonal vector being $x_j.$

Note that for any tuple $(x_1,\ldots, x_d)$ of vectors in $\mathbb{C}^d,$ the tuple of Varopoulos operators  $\mathbf{V}=(V_1,\ldots, V_d)$ of order $k$ is a commuting tuple. 
Moreover, $V_j$ is a contraction if and only if $\|x_j\|_2\le 1.$ The tuple 
 $\mathbf{V}$ is both a row and a spherical
contraction if and only if $\sum_{j=1}^d\|x_j\|_2^2\leq1.$
Let $\mathcal{V}_{d,k}$ denote the set of all commuting contractive $d$-tuple of Varopoulos operators of order $k.$ Let $\mathfrak{R}_{d,k}$ denote the set of all row contractive $d$-tuple of Varopoulos operators of order $k.$ 
\end{definition}
 \begin{theorem}\label{Varopoulos-Constant-order-k} 
Let $p(z_1,\ldots,z_d)$ be a homogeneous polynomial of degree $k$, and let $V_1,\ldots,V_d$ be contractive Varopoulos operators of order $k$. Then 
$$\| p(V_1,\ldots,V_d) \|
  \leq K^{\mathbb{C}}_G \frac{k^{k/2} (k+1)^{(k+1)/2}}{k!2^{\frac{k+2}{2}}}\,\|p\|_{\infty,\mathbb{D}^d}.$$ 
  \end{theorem} 
  \begin{proof}
  The proof follows from \cite[A(r) Inequality]{TG} put together with Theorem \ref{Injective-tensor-norm}.    
  \end{proof}
 Note that the bound in Theorem \ref{Varopoulos-Constant-order-k} is independent of the number of variables $d.$ Therefore the word `homogeneous' can be removed from the statement of the same theorem. In the following remark we further note that the ratio $\frac{\|p(T_1,\dots,T_d)\|}{\|p(V_1,\dots,V_d)\|}$ goes to infinity if we choose to vary the permissible tuples $(T_1,\ldots T_d),$ $(V_1,\ldots, V_d),$ and the polynomial $p$ of degree $k$ for any $k\ge 3.$
\begin{remark}
Note that the ratio $\frac{\|p(T_1,\dots,T_d)\|}{\|p(V_1,\dots,V_d)\|}$ is at least $1$ if we vary $\boldsymbol T \in \mathscr{C}_d,$ $\boldsymbol V\in\mathcal{V}_{d,k}$ and $p\in\mathcal P_k.$
Take $k=3$ (as the arguments repeats for the general $k\ge 3$). Suppose 
\begin{eqnarray}\label{if-ratio-is-bounded}
\frac{\sup_{\boldsymbol{T}\in \mathscr{C}_d}\|p(T_1,\dots,T_d)\|}{\sup_{\boldsymbol V\in\mathcal{V}_{d,k}}\|p(V_1,\dots,V_d)\|}\leq M \quad \mbox{ for all }p\in \mathcal P_k \mbox{ and } d\in\mathbb N,
\end{eqnarray}
for some positive real number $M.$ Then using Theorem \ref{Varopoulos-Constant-order-k}, we can say that 
\begin{eqnarray}\label{False-assumption-bound-General}
    \sup_{\boldsymbol{T}\in \mathscr C_d}\|p(T_1,\dots,T_d)\|\leq M K^{\mathbb{C}}_G \frac{k^{k/2} (k+1)^{(k+1)/2}}{k!2^{\frac{k+2}{2}}} \|p\|_\infty, \quad \mbox{ for any }p\in \mathcal P_k \mbox{ and } d\in\mathbb N.
\end{eqnarray}
On the other hand according to \cite[Theorem 1]{V1}, for any constant $K>0$ there exists a $d\in\mathbb N,$ polynomial $p$ of degree $3$ in $d$ variables and a tuple $\boldsymbol T\in \mathscr{C}_d$ such that $\|p(\boldsymbol T)\|> K \|p\|_\infty.$ This immediately contradicts \eqref{False-assumption-bound-General}. Hence \eqref{if-ratio-is-bounded} is not possible for any real number $M.$
\end{remark}
Next natural question is: what happens to the ratio $\frac{\sup_{\boldsymbol{T}\in \mathscr{C}_d}\|P(T_1,\dots,T_d)\|}{\sup_{\boldsymbol V\in\mathcal V_{d,k}}\|P(V_1,\dots,V_d)\|}$ for a fixed $d?$ Addressing this question is the content of the following remark.
\begin{remark}
Theorem 1.1(b) of \cite{P.G Dixon} states that for any fixed $d\in\mathbb N,$ there exists a homogeneous polynomial of degree $3$ such that $\sup\left\{\|p(\mathbf{T})\| :\mathbf{T}\in \mathscr{C}_d \right\}> \frac{1}{11}d^{1/2}\|p\|_{\infty}.$ 
From Theorem~\ref{Varopoulos-Constant-order-k}, it follows that
\[
\frac{\sup_{\boldsymbol T\in\mathscr C_d} \|p(T_1,\dots,T_d)\|}{\sup_{\boldsymbol V\in \mathcal{V}_{d,3}} \|p(V_1,\dots,V_d)\|}
\;\geq\; \frac{1}{11\sqrt{6}\,K^{\mathbb{C}}_G}\, d^{1/2}.
\]
\end{remark}

\section{Ball version of Grothendieck inequality}
The ball version of Grothendieck's inequality, as discussed by Mantero and Tonge in
\cite[Lemma 4.1]{ManteroTonge1976}, plays a crucial role in
establishing von Neumann inequalities for homogeneous polynomials
of degree two in an arbitrary number of variables. Building on
this result, we further develop a multidimensional version of
Grothendieck's inequality using symmetric product (Definition \ref{Symmetric-product}). This generalized Grothendieck's inequality enables
us to establish a ball version of von Neumann's inequality for
the class of row contractive Varopoulos operators. The following
theorem presents this result.

\begin{theorem}\cite[Lemma 4.1]{ManteroTonge1976}\label{ball-Grothendieck-inequality}
Let $1\leq p\leq 2$, and suppose that $x_i$ ($1\leq i\leq N$)
and $y_m$ ($1\leq j\leq M$) are elements of a Hilbert space satisfying
\begin{eqnarray}\label{contraction-condition}
\sum_i \|x_i\|^p\leq 1
\qquad\text{and}\qquad
\sum_j \|y_j\|^p\leq 1.
\end{eqnarray}
Then, for any matrix $(a_{ij})$, we have
\[
\left|\sum_i\sum_j a_{ij}\langle x_i,y_j\rangle\right|
\leq
\sup_{s,t}
\left|\sum_i\sum_j a_{ij}s_it_j\right|,
\]
where the supremum is taken over all
$s=(s_i)$ and $t=(t_j)$ satisfying
\[
\left(\sum_i |s_i|^p\right)^{1/p}\leq 1
\qquad\text{and}\qquad
\left(\sum_j |t_j|^p\right)^{1/p}\leq 1.
\]
\end{theorem}

The following theorem provides a generalization of the Theorem \ref{ball-Grothendieck-inequality}, which is stated as follows.

\begin{theorem}\label{Ball-version-of-multidimentional-Grothendieck-inequality}
Let $1 \leq p \leq 2$ and let $k \geq 2$ be an integer. Let $(a_{i_1 i_2 \dots i_k})_{1 \leq i_1, \dots, i_k \le d}$ be a  $k$-tensor. For each $m \in \{1, \dots, k\}$, let $\{x_{i}^{(m)}\}_{i=1}^d \subset \ell^2$ be a collection of vectors satisfying
\[
\sum_{i=1}^d \|x_{i}^{(m)}\|_2^p \leq 1 \quad \text{for each } m \in \{1, \dots, k\}.
\]
Then
\[
\left| \sum_{i_1, \ldots, i_k = 1}^d a_{i_1 \dots i_k} [x_{i_1}^{(1)} , x_{i_2}^{(2)}, \dots, x_{i_k}^{(k)}] \right| \leq \sup \left| \sum_{i_1, \dots, i_k = 1}^d a_{i_1 \dots i_k} s_{i_1}^{(1)} s_{i_2}^{(2)} \cdots s_{i_k}^{(k)} \right|,
\]
where the supremum is taken over all scalar sequences $s_m = (s_{1}^{(m)}, \dots, s_{d}^{(m)}) \in \mathbb{C}^d$ satisfying $\left(\sum_{i=1}^d |s_{i}^{(m)}|^p\right)^{1/p} \le 1$ for all $m \in \{1, \dots, k\}$.
\end{theorem}

\begin{proof}
Define the scalar multilinear bound as follows
\[
M = \sup \left\{ \left| \sum_{i_1, \dots, i_k = 1}^d a_{i_1 \dots i_k} \prod_{m=1}^k s_{i_m}^{(m)} \right| : \sum_{i_{m}=1}^d |s_{i_{m}}^{(m)}|^p \leq 1 \text{ for } 1 \leq m \leq k \right\}.
\]
By multilinear scaling, for arbitrary scalar sequences $s_{i_{m}}^{(m)}$, we have
\begin{eqnarray}\label{estimation-1}
\left| \sum_{i_1, \dots, i_k = 1}^d a_{i_1 \dots i_k} \prod_{m=1}^k s_{i_m}^{(m)} \right| \leq M \prod_{m=1}^k \left( \sum_{i_{m}=1}^d |s_{i_{m}}^{(m)}|^p \right)^{1/p}.
\end{eqnarray}

Expanding the product $[x_{i_1}^{(1)}, \dots, x_{i_k}^{(k)}] = \sum_{l=1}^{\infty} \prod_{m=1}^k x_{i_m}^{(m)}(l)$ and interchanging the order of summation gives
\[
\left| \sum_{i_1, \dots, i_k = 1}^d a_{i_1 \dots i_k} [x_{i_1}^{(1)}, \dots, x_{i_k}^{(k)}] \right| \leq \sum_{l=1}^{\infty} \left| \sum_{i_1, \dots, i_k = 1}^d a_{i_1 \dots i_k} \prod_{m=1}^k x_{i_m}^{(m)}(l) \right|.
\]

Applying \eqref{estimation-1} with $s_{i_{m}}^{(m)} = x_{i_m}^{(m)}(l)$ yields
\begin{eqnarray}\label{Ball-Grothendieck}
\left| \sum_{i_1, \dots, i_k = 1}^d a_{i_1 \dots i_k} [x_{i_1}^{(1)}, \dots, x_{i_k}^{(k)}] \right| \leq M \sum_{l=1}^{\infty} \prod_{m=1}^k \left( \sum_{i_{m}=1}^d |x_{i_m}^{(m)}(l)|^p \right)^{1/p}.
\end{eqnarray}

By Generalized H\"older's inequality  with $k$ equal conjugate exponents $r_1 = r_2 = \dots = r_k = k$ (satisfying $\sum_{m=1}^k \frac{1}{k} = 1$), we obtain
\[
\sum_{l=1}^{\infty} \prod_{m=1}^k \left( \sum_{i_m=1}^d |x_{i_m}^{(m)}(l)|^p \right)^{1/p} \leq \prod_{m=1}^k \left[ \sum_{l=1}^{\infty} \left( \sum_{i_{m}=1}^d |x_{i_m}^{(m)}(l)|^p \right)^{k/p} \right]^{1/k}.
\]

Since $1 \leq p \leq 2$ and $k \geq 2$, we have $\frac{2}{p} \geq 1$. For any sequence of non-negative terms $c_l = \sum_{i_{m}=1}^d |x_{i_m}^{(m)}(l)|^p$, the monotonicity of $\ell_q$ norms ($\ell_{2/p} \hookrightarrow \ell_{k/p}$ for $k \geq 2$) gives
\[
\left[ \sum_{l=1}^{\infty} c_l^{k/p} \right]^{1/k} \leq \left[ \sum_{l=1}^{\infty} c_l^{2/p} \right]^{1/2} = \left[ \sum_{l=1}^{\infty} \left( \sum_{i_{m}=1}^d |x_{i_m}^{(m)}(l)|^p \right)^{2/p} \right]^{1/2}.
\]

Applying Minkowski's inequality for $\ell_{2/p}$ (since $\frac{2}{p} \geq 1$) to exchange the sum over $d$ and the sum over $i$ yields
\begin{eqnarray}\label{estimation}
\left[ \sum_{l=1}^{\infty} \left( \sum_{i_{m}=1}^d |x_{i_m}^{(m)}(l)|^p \right)^{2/p} \right]^{1/2} \leq \left[ \sum_{i_{m}=1}^d \left( \sum_{l=1}^{\infty} |x_{i_m}^{(m)}(l)|^2 \right)^{p/2} \right]^{1/p} = \left( \sum_{i=1}^d \|x_{i}^{(m)}\|_2^p \right)^{1/p} \le 1.
\end{eqnarray}

Substituting  bound \eqref{estimation}  for all $m \in \{1, \dots, k\}$ in \eqref{Ball-Grothendieck} , we have 
\[
\left| \sum_{i_1, \dots, i_k = 1}^d a_{i_1 \dots i_k} [x_{i_1}^{(1)}, \dots, x_{i_k}^{(k)}] \right| \leq M.
\]
This completes the proof.
\end{proof}
\begin{theorem}\cite[Proposition 3]{DeyGuptaKumar2026}
 Let \( \boldsymbol T = (T_1, \ldots, T_d) \) be a doubly contractive operator tuple on a Hilbert space $\mathcal{H}$ and $p\in \mathbb{C}[z_1,\ldots ,z_d]$  be a homogeneous polynomial of degree $2$ in $d$-variables then $\|p(T)\|\leq \|p\|_{\mathbb{B}_d,\infty}$. 
 \end{theorem} 
\begin{proof}
An alternative proof of the theorem follows from an application of Theorem
\ref{ball-Grothendieck-inequality}. Let $x,y\in\mathcal H$ be arbitrary
unit vectors. Then
\[
\left|\langle p(T_1,\ldots,T_d)x,y\rangle\right|
=
\left|\sum_{i,j=1}^d a_{ij}
\langle T_i x,T_j^*y\rangle\right|.
\]
Since $(T_1,\ldots,T_d)$ is both a row contraction and a column
contraction, we have
\[
\sum_{i=1}^d\|T_i x\|^2\leq 1
\qquad\text{and}\qquad
\sum_{j=1}^d\|T_j^*y\|^2\leq 1.
\]
Consequently, the contraction condition \eqref{contraction-condition}
is satisfied for $p=2$. Thus, an application of Theorem
\ref{ball-Grothendieck-inequality} yields the desired estimate.

Furthermore, by \cite[Proposition 1.43,1.44]{Di},
applied to the associated bilinear form, we obtain $\sup_{s,t}
\left|\sum_{i,j}a_{ij}s_it_j\right|
=
\sup_t
\left|\sum_{i,j}a_{ij}t_it_j\right|,$
where the first supremum is taken over all sequences
$s=(s_i)$ and $t=(t_j)$ satisfying $\left(\sum_i|s_i|^2\right)^{1/2}\leq 1$ and $\left(\sum_j|t_j|^2\right)^{1/2}\leq 1,$
respectively, while the second supremum is taken over all sequences
$t=(t_i)$ satisfying $\left(\sum_i|t_i|^2\right)^{1/2}\leq 1.$
This completes the proof.
\end{proof}
\begin{theorem}\label{ Ball-version-of-Varopoulos-Constant-order-k} 
Let $p(z_1,\ldots,z_d)$ be a homogeneous polynomial of degree $k$, and let $V_1,\ldots,V_d$ be row contractive Varopoulos operators of order $k$. Then 
$$\| p(V_1,\ldots,V_d) \|
  \leq \|p\|_{\mathbb{B}_d ,\infty}.$$ 
  \end{theorem} 
\begin{proof}
By the definition of $p(V_1,\ldots,V_d)$, we have
\[
\|p(V_1,\ldots,V_d)\|
=
\sup_{\substack{\sum_{j=1}^d\|x_j\|_{2}^2\leq 1\\
}}
\left|
\sum_{i_1,\ldots,i_k=1}^d
a_{i_1\ldots i_k}
[x_{i_1},\ldots,x_{i_k}]
\right|.
\]
We now apply Theorem \ref{Ball-version-of-multidimentional-Grothendieck-inequality},
together with the following identity:
\[
\begin{aligned}
&\sup_{\substack{s_m=(s_{1}^{(m)},\ldots,s_{d}^{(m)})\in\mathbb C^d,\\
\left(\sum_{i=1}^d|s_{i}^{(m)}|^2\right)^{1/2}\leq 1,\;
m=1,\ldots,k}}
\left|
\sum_{i_1,\ldots,i_k=1}^d
a_{i_1\ldots i_k}
s_{i_1}^{(1)}\cdots s_{i_k}^{(k)}
\right|\\
&=
\sup_{\substack{(s_{1},\ldots,s_{d})\in\mathbb C^d,\\
\left(\sum_{i=1}^d|s_{i}|^2\right)^{1/2}\leq 1}}
\left|
\sum_{i_1,\ldots,i_k=1}^d
a_{i_1\ldots i_k}
s_{i_1}\cdots s_{i_k}
\right|.
\end{aligned}
\]
This identity follows from \cite[Propositions 1.43 and 1.44]{Di}. This completes the proof.
\end{proof}

\subsection*{Acknowledgments} 
I would like to sincerely thank my supervisor Dr. Rajeev Gupta for his invaluable guidance and constant support during the preparation of this article.

\bibliographystyle{amsplain}

\end{document}